# EQUIVALENCE THEORY FOR DENSITY ESTIMATION, POISSON PROCESSES AND GAUSSIAN WHITE NOISE WITH DRIFT

By Lawrence D. Brown[1], Andrew V. Carter, Mark G. Low[2]
and Cun-Hui Zhang[3]

*University of Pennsylvania, University of California, Santa Barbara,
University of Pennsylvania and Rutgers University*

This paper establishes the global asymptotic equivalence between a Poisson process with variable intensity and white noise with drift under sharp smoothness conditions on the unknown function. This equivalence is also extended to density estimation models by Poissonization. The asymptotic equivalences are established by constructing explicit equivalence mappings. The impact of such asymptotic equivalence results is that an investigation in one of these nonparametric models automatically yields asymptotically analogous results in the other models.

**1. Introduction.** The purpose of this paper is to give an explicit construction of global asymptotic equivalence in the sense of Le Cam (1964) between a Poisson process with variable intensity and white noise with drift. The construction is extended to density estimation models. It yields asymptotic solutions to both density estimation and Poisson process problems based on asymptotic solutions to white noise with drift problems and vice versa.

*Density estimation model.* A random vector $\mathbf{V}_n^\star$ of length $n$ is observed such that $\mathbf{V}_n^\star \equiv (V_1^\star, \ldots, V_n^\star)$ is a sequence of i.i.d. variables with a common density $f \in \mathcal{F}$.

Received April 2002; revised February 2004.
[1]Supported by NSF Grant DMS-99-71751.
[2]Supported by NSF Grant DMS-03-06576.
[3]Supported by NSF Grants DMS-01-02529 and DMS-02-03086.

*Poisson process.* A random vector of random length $\{N, \mathbf{X}_N\}$ is observed such that $N \equiv N_n$ is a Poisson variable with $EN = n$ and that given $N = m$, $\mathbf{X}_N = \mathbf{X}_m \equiv (X_1, \ldots, X_m)$ is a sequence of i.i.d. variables with a common density $f \in \mathcal{F}$. The resulting observations are then distributed as a Poisson process with intensity function $nf$.

*White noise.* A Gaussian process $Z^\star \equiv Z_n^\star \equiv \{Z_n^\star(t), 0 \leq t \leq 1\}$ is observed such that

$$(1.1) \qquad Z_n^\star(t) \equiv \int_0^t \sqrt{f(x)}\, dx + \frac{B^\star(t)}{2\sqrt{n}}, \qquad 0 \leq t \leq 1,$$

with a standard Brownian motion $B^\star(t)$ and an unknown probability density function $f \in \mathcal{F}$ in $[0, 1]$.

*Asymptotic equivalence.* For any two experiments $\xi_1$ and $\xi_2$ with a common parameter space $\Theta$, $\Delta(\xi_1, \xi_2; \Theta)$ denotes Le Cam's distance [cf., e.g., Le Cam (1986) or Le Cam and Yang (1990)] defined as

$$\Delta(\xi_1, \xi_2; \Theta) \equiv \sup_L \max_{j=1,2} \sup_{\delta^{(j)}} \inf_{\delta^{(k)}} \sup_{\theta \in \Theta} |E_\theta^{(j)} L(\theta, \delta^{(j)}) - E_\theta^{(k)} L(\theta, \delta^{(k)})|,$$

where (a) the first supremum is taken over all decision problems with loss function $\|L\|_\infty \leq 1$, (b) given the decision problem and $j = 1, 2$, $k \equiv 3 - j$ ($k = 2$ for $j = 1$ and $k = 1$ for $j = 2$) the "maximin" value of the maximum difference in risks over $\Theta$ is computed over all (randomized) statistical procedures $\delta^{(\ell)}$ for $\xi_\ell$ and (c) the expectations $E_\theta^{(\ell)}$ are evaluated in experiments $\xi_\ell$ with parameter $\theta$, $\ell = j, k$. The statistical interpretation of the Le Cam distance is as follows: If $\Delta(\xi_1, \xi_2; \Theta) < \varepsilon$, then for any decision problem with $\|L\|_\infty \leq 1$ and any statistical procedure $\delta^{(j)}$ with the experiment $\xi_j$, $j = 1, 2$, there exists a (randomized) procedure $\delta^{(k)}$ with $\xi_k$, $k = 3 - j$, such that the risk of $\delta^{(k)}$ evaluated in $\xi_k$ nearly matches (within $\varepsilon$) that of $\delta^{(j)}$ evaluated in $\xi_j$.

Two sequences of experiments $\{\xi_{1,n}, n \geq 1\}$ and $\{\xi_{2,n}, n \geq 1\}$, with a common parameter space $\mathcal{F}$, are asymptotically equivalent if

$$\Delta(\xi_{1,n}, \xi_{2,n}; \mathcal{F}) \to 0 \qquad \text{as } n \to \infty.$$

The interpretation is that the risks of corresponding procedures converge.

A key result of Le Cam (1964) is that this equivalence of experiments can be characterized using random transformations between the probability spaces. A random transformation, $T(X, U)$ which maps observations $X$ into the space of observations $Y$ (with possible dependence on an independent, uninformative random component $U$) also maps distributions in $\xi_1$ to approximations of the distributions in $\xi_2$ via $\mathbf{P}_\theta^{(1)} T \approx \mathbf{P}_\theta^{(2)}$. For the mapping between the Poisson and Gaussian processes we shall restrict ourselves to transformations $T$ with deterministic inverses, $T^{-1}(T(X,U)) = X$. The experiments are asymptotically equivalent if the total-variation distance between $\mathbf{P}_\theta^{(2)}$ and the distribution of $T$ under $\mathbf{P}_\theta^{(1)}$ converges to 0 uniformly in $\theta$. As explained in Brown and Low (1996) and Brown, Cai, Low and Zhang (2002), knowing an appropriate $T$ allows explicit construction of estimation procedures in $\xi_1$ by applying statistical procedures from $\xi_2$ to $T(X,U)$.

In general, asymptotic equivalence also implies a transformation from the $\mathbf{P}_\theta^{(2)}$ to the $\mathbf{P}_\theta^{(1)}$ and the corresponding total-variation distance bound. However, in the case of the equivalence between the Poisson process and white noise with drift, by requiring that the transformation be invertible, we have saved ourselves a step. The transformation in the other direction is $T^{-1}$, and

$$\|\mathbf{P}_\theta^{(1)}T - \mathbf{P}_\theta^{(2)}\| \geq \|\mathbf{P}_\theta^{(1)}TT^{-1} - \mathbf{P}_\theta^{(2)}T^{-1}\| = \|\mathbf{P}_\theta^{(1)} - \mathbf{P}_\theta^{(2)}T^{-1}\|.$$

Therefore, it is sufficient if $\sup_\theta \|\mathbf{P}_\theta^{(1)}T - \mathbf{P}_\theta^{(2)}\| \to 0$.

The equivalence mappings $T_n$ constructed in this paper from the sample space of the Poisson process to the sample space of the white noise are invertible randomized mappings such that

(1.2) $$\sup_{f \in \mathcal{F}} H_f(T_n(N, \mathbf{X}_N), Z_n^\star) \to 0$$

under certain conditions on the family $\mathcal{F}$. Here $H_f(Z_1, Z_2)$ denotes the Hellinger distance of stochastic processes or random vectors $Z_1$ and $Z_2$ living in the same sample space, when the true unknown density is $f$. Since $T_n$ are invertible randomized mappings, $T_n(N, \mathbf{X}_N)$ are sufficient statistics for the Poisson processes and their inverses $T_n^{-1}$ are necessarily many-to-one deterministic mappings. Similar considerations apply for the mapping of the density estimation problem to the white noise with drift problem, although in that case there are two mappings, one from the density estimation to the white noise with drift model and another from the white noise with drift model back to the density estimation model. These mappings are given in Section 2.

There have recently been several papers on the global asymptotic equivalence of nonparametric experiments. Brown and Low (1996) established global asymptotic equivalence of the white noise problem with unknown drift $f$ to a nonparametric regression problem with deterministic design and unknown regression $f$ when $f$ belongs to a Lipschitz class with smoothness index $\alpha > \frac{1}{2}$. It has also been demonstrated that such nonparametric problems are typically asymptotically nonequivalent when the unknown $f$ belongs to larger classes, for example, with smoothness index $\alpha \leq \frac{1}{2}$. Brown and Low (1996) showed the asymptotic nonequivalence between the white noise problem and nonparametric regression with deterministic design for $\alpha \leq \frac{1}{2}$, Efromovich and Samarov (1996) showed that the asymptotic equivalence may fail when $\alpha < \frac{1}{4}$. Brown and Zhang (1998) showed the asymptotic nonequivalence for $\alpha \leq \frac{1}{2}$ between any pair of the following four experiments: white noise, density problem, nonparametric regression with random design, and nonparametric regression with deterministic design. In Brown, Cai, Low and Zhang (2002) the asymptotic equivalence for nonparametric regression with random design was shown under Besov constraints which include Lipschitz classes with any smoothness index $\alpha > \frac{1}{2}$. Gramma and Nussbaum (1998) solved the fixed-design nonparametric regression problem for nonnormal errors. Milstein and Nussbaum (1998) showed that some diffusion problems can be approximated by discrete versions that are nonparametric autoregression models, and Golubev and Nussbaum (1998) established a discrete Gaussian approximation to the problem of estimating the spectral density of a stationary process.

Most closely related to this paper is the work in Nussbaum (1996) where global asymptotic equivalence of the white noise problem to the nonparametric density problem with unknown density $g = f^2/4$ is shown. In this paper the global asymptotic equivalence was established under the following smoothness assumption: $f$ belongs to the Lipschitz classes with smoothness index $\alpha > \frac{1}{2}$.

*The parameter spaces.* The class of functions $\mathcal{F}$ will be assumed throughout to be densities with respect to Lebesgue measure on $[0,1]$ that are uniformly bounded away from 0. The smoothness conditions on $\mathcal{F}$ can be described in terms of Haar basis functions of the densities. Let

$$(1.3) \quad \theta_{k,\ell} \equiv \theta_{k,\ell}(f) \equiv \int f \phi_{k,\ell}, \qquad \ell = 0, \ldots, 2^k - 1, \ k = 0, 1, \ldots,$$

be the Haar coefficients of $f$, where

$$(1.4) \qquad \phi_{k,\ell} \equiv 2^{k/2}(\mathbb{1}_{I_{k+1,2\ell}} - \mathbb{1}_{I_{k+1,2\ell+1}})$$

are the Haar basis functions with $I_{k,\ell} \equiv [\ell/2^k, (\ell+1)/2^k)$. The convergence of the Hellinger distance in (1.2) is established via an inequality in Theorem 3 in terms of the tails of the Besov norms $\|f\|_{1/2,2,2}$ and $\|f\|_{1/2,4,4}$ of the Haar coefficients $\theta_{k,\ell} \equiv \theta_{k,\ell}(f)$ in (1.3).

The Besov norms $\|f\|_{\alpha,p,q}$ for the Haar coefficients, with smoothness index $\alpha$ and shape parameters $p$ and $q$, are defined by

$$(1.5) \quad \|f\|_{\alpha,p,q} \equiv \left[ \left| \int_0^1 f \right|^q + \sum_{k=0}^{\infty} \left\{ 2^{k(\alpha+1/2-1/p)} \left( \sum_{\ell=0}^{2^k-1} |\theta_{k,\ell}(f)|^p \right)^{1/p} \right\}^q \right]^{1/q}.$$

Let $\bar{f}_k$ be the piecewise average of $f$ at resolution level $k$, that is, the piecewise constant function defined by

$$\bar{f}_k \equiv \bar{f}_k(t) \equiv \sum_{\ell=0}^{2^k-1} \mathbb{1}\{t \in I_{k,\ell}\} 2^k \int_{I_{k,\ell}} f. \tag{1.6}$$

Since $\|\bar{f}_k - \bar{f}_{k+1}\|_p^p = \int |\sum_\ell \theta_{k,\ell} \phi_{k,\ell}|^p = \sum_\ell |\theta_{k,\ell}|^p 2^{k(p/2-1)}$, (1.5) can be written as

$$\|f\|_{\alpha,p,q} \equiv \left\{ |\bar{f}_0|^q + \sum_{k=0}^{\infty} (2^{k\alpha} \|\bar{f}_k - \bar{f}_{k+1}\|_p)^q \right\}^{1/q},$$

and its tail at resolution level $k_0 \geq 0$ is $\|f - \bar{f}_{k_0}\|_{\alpha,p,q}$, $k_0 \geq 0$, with

$$\|f - \bar{f}_{k_0}\|_{\alpha,p,q}^q = \sum_{k=k_0}^{\infty} \left\{ 2^{k(\alpha+1/2-1/p)} \left( \sum_{\ell=1}^{2^k-1} |\theta_{k,\ell}|^p \right)^{1/p} \right\}^q. \tag{1.7}$$

Let $B(\alpha, p, q)$ be the Besov space

$$B(\alpha, p, q) = \{f : \|f\|_{\alpha,p,q} < \infty\}.$$

The following two theorems on the equivalence of white noise with drift, density estimation and Poisson estimation models are corollaries of our main result, Theorem 3, which bounds the squared Hellinger distance between particular invertible randomized mappings of the Poisson process and white noise with drift models. The randomized mappings are given in Section 2. Proofs of these theorems are given in the Appendix.

THEOREM 1. *Let $Z_n^\star$, $\{N, \mathbf{X}_N\}$ and $\mathbf{V}_n^\star$ be the Gaussian process, Poisson process and density estimation experiments, respectively. Suppose that $\mathcal{H}$ is compact in both $B(1/2, 2, 2)$ and $B(1/2, 4, 4)$ and that $\mathcal{H} \subseteq \{f : \inf_{0 < x < 1} f(x) \geq \varepsilon_0\}$ for some $\varepsilon_0 > 0$. Then*

$$\lim_{n \to \infty} \Delta(Z_n^\star, \{N, \mathbf{X}_N\}; \mathcal{H}) = 0 \tag{1.8}$$

*and*

$$\lim_{n \to \infty} \Delta(Z_n^\star, \mathbf{V}_n^\star; \mathcal{H}) = 0. \tag{1.9}$$

Our construction also shows that asymptotic equivalence holds for a class $\mathcal{F}$ if $\mathcal{F}$ is bounded in the Lipschitz norm with smoothness index $\beta$ and compact in the Sobolev norm with smoothness index $\alpha \geq \beta$ such that $\alpha + \beta \geq 1$, $\alpha \geq \frac{3}{4}$ or $\beta > \frac{1}{2}$.

For $0 < \beta \leq 1$ the Lipschitz norm $\|f\|_\beta^{(L)}$ and Sobolev norm $\|f\|_\alpha^{(S)}$ are defined by

$$\|f\|_\beta^{(L)} \equiv \sup_{0 \leq x < y \leq 1} \frac{|f(x) - f(y)|}{|x-y|^\beta}, \qquad \|f\|_\alpha^{(S)} \equiv \sum_{n=-\infty}^{\infty} n^{2\alpha} |c_n(f)|^2,$$

where $c_n(f) \equiv \int_0^1 f(x) e^{-in2\pi x} \, dx$ are the Fourier coefficients of $f$.

THEOREM 2. *Let $Z_n^\star$, $\{N, \mathbf{X}_N\}$ and $\mathbf{V}_n^\star$ be the Gaussian process, Poisson process and density estimation experiments, respectively, and let $\mathcal{F}$ be bounded in the Lipschitz norm with smoothness index $\beta$ and compact in the Sobolev norm with smoothness index $\alpha \geq \beta$. Suppose $\mathcal{F} \subseteq \{f : \inf_{0 < x < 1} f(x) \geq \varepsilon_0\}$ for some $\varepsilon_0 > 0$. Then if $\alpha + \beta \geq 1$, $\alpha \geq \frac{3}{4}$ or $\beta > \frac{1}{2}$,*

$$\lim_{n \to \infty} \Delta(Z_n^\star, \{N, \mathbf{X}_N\}; \mathcal{F}) = 0$$

*and*

$$\lim_{n \to \infty} \Delta(Z_n^\star, \mathbf{V}_n^\star; \mathcal{F}) = 0.$$

**2. The equivalence mappings.** This section describes in detail the mappings which provide the asymptotic equivalence claimed in this paper. The fact that these mappings yield asymptotic equivalence is a consequence of our major result, Theorem 3. The construction is broken into several stages.

From observations of the white noise (1.1), define random vectors

(2.1) $\quad \overline{\mathbf{Z}}_k^\star \equiv \{\overline{Z}_{k,\ell}^\star, 0 \leq \ell < 2^k\}, \qquad \overline{Z}_{k,\ell}^\star \equiv 2^k \left\{ Z^\star \left( \frac{\ell+1}{2^k} \right) - Z^\star \left( \frac{\ell}{2^k} \right) \right\},$

(2.2) $\quad \mathbf{W}_k^\star \equiv \{W_{k,\ell}^\star, 0 \leq \ell < 2^k\}, \qquad W_{k,2\ell}^\star \equiv -W_{k,2\ell+1}^\star \equiv (\overline{Z}_{k,2\ell}^\star - \overline{Z}_{k,2\ell+1}^\star)/2.$

Let $k_0 \equiv k_{0,n}$ be suitable integers with $\lim_{n \to \infty} k_{0,n} = \infty$. Following Brown, Cai, Low and Zhang (2002), we construct equivalence mappings by finding the counterparts of $\overline{\mathbf{Z}}_{k_0}^\star$ and $\mathbf{W}_k^\star$, $k > k_0$, with the Poisson process $(N, \mathbf{X}_N)$, to strongly approximate the Gaussian variables.

It can be easily verified from (1.1) that $\{\overline{Z}_{k_0,\ell}^\star, 0 \leq \ell < 2^{k_0}, W_{k,2\ell}^\star, 0 \leq \ell < 2^{k-1}, k > k_0\}$ are uncorrelated normal random variables with

(2.3)
$$E\overline{Z}_{k,\ell}^\star = h_{k,\ell} \equiv 2^k \int_{I_{k,\ell}} h, \qquad h \equiv \sqrt{f},$$
$$\sqrt{\operatorname{Var}(\overline{Z}_{k,\ell}^\star)} = \sigma_k \equiv \sqrt{2^k/(4n)},$$

for $\ell = 0, \ldots, 2^k - 1$, and for $\ell = 0, \ldots, 2^{k-1} - 1$,

(2.4)
$$EW_{k,2\ell}^\star = \tfrac{1}{2}(h_{k,2\ell} - h_{k,2\ell+1}) = \sqrt{2^{k-1}} \int h \phi_{k-1,\ell},$$
$$\sqrt{\operatorname{Var}(W_{k,2\ell}^\star)} = \sigma_{k-1}.$$

Let $\widetilde{\mathbf{U}} = \{\widetilde{U}_{k,\ell}, k \geq k_0, \ell \geq 0\}$ be a sequence of i.i.d. uniform variables in $[-1/2, 1/2)$ independent of $(N, \mathbf{X}_N)$. For $k = 0, 1, \ldots$ and $\ell = 0, \ldots, 2^k - 1$ define

$$(2.5) \qquad \mathbf{N}_k \equiv \{N_{k,\ell}, 0 \leq \ell < 2^k\}, \qquad N_{k,\ell} \equiv \#\{X_i : X_i \in I_{k,\ell}\}.$$

We shall approximate $\overline{\mathbf{Z}}_k^\star$ in (2.1) in distribution by

$$(2.6) \qquad \begin{aligned} \overline{\mathbf{Z}}_k &\equiv \{\overline{Z}_{k,\ell}, 0 \leq \ell < 2^k\}, \\ \overline{Z}_{k,\ell} &\equiv 2\sigma_k \operatorname{sgn}(N_{k,\ell} + \widetilde{U}_{k,\ell}) \sqrt{|N_{k,\ell} + \widetilde{U}_{k,\ell}|}, \end{aligned}$$

at the initial resolution level $k = k_0$. Since $N_{k,\ell}$ are Poisson variables with

$$(2.7) \qquad \lambda_{k,\ell} \equiv E N_{k,\ell} = \frac{n}{2^k} f_{k,\ell} = \frac{f_{k,\ell}}{4\sigma_k^2}, \qquad f_{k,\ell} \equiv 2^k \int_{I_{k,\ell}} f,$$

by Taylor expansion and central limit theory

$$\overline{Z}_{k,\ell} \approx 2\sigma_k \left( \sqrt{\lambda_{k,\ell}} + \frac{N_{k,\ell} - \lambda_{k,\ell}}{2\lambda_{k,\ell}^{1/2}} \right) \approx N(\sqrt{f_{k,\ell}}, \sigma_k^2)$$

as $\lambda_{k,\ell} \to \infty$, compared with (2.3). Note that $\sqrt{f_{k,\ell}} \approx h_{k,\ell}$ under suitable smoothness conditions on $f$, in view of (2.3) and (2.7). The Poisson variables $N_{k,\ell}$ can be fully recovered from $\overline{Z}_{k,\ell}$, while the randomization turns $N_{k,\ell}$ into continuous variables.

Approximation of $W_{k,\ell}^\star$ for $k > k_0$ is more delicate, since the central limit theorem is not sufficiently accurate at high resolution levels. Let $F_m$ be the cumulative distribution function of the independent sum of a binomial variable $\widetilde{X}_{m,1/2}$ with parameter $(m, \frac{1}{2})$ and a uniform variable $\widetilde{U}$ in $[-\frac{1}{2}, \frac{1}{2})$,

$$(2.8) \qquad F_m(x) \equiv P\{\widetilde{X}_{m,1/2} + \widetilde{U} \leq x\},$$

with $F_0$ being the uniform distribution in $[-\frac{1}{2}, \frac{1}{2})$. Let $\Phi$ be the $N(0,1)$ cumulative distribution. We shall approximate $\mathbf{W}_k^\star$ by using a quantile transformation of randomized versions of the Poisson random variables. More specifically, let

$$(2.9) \quad \mathbf{W}_k \equiv \{W_{k,\ell}, 0 \leq \ell < 2^k\}, \qquad W_{k,2\ell} \equiv \sigma_{k-1} \Phi^{-1}(F_{N_{k-1,\ell}}(N_{k,2\ell} + \widetilde{U}_{k,2\ell}))$$

with $W_{k,2\ell} \equiv -W_{k,2\ell+1}$, $\ell = 0, \ldots, 2^{k-1} - 1$, and the $\sigma_k$ in (2.3). Given $N_{k-1,\ell} = m$,

$$(2.10) \quad N_{k,2\ell} \sim \operatorname{Bin}(m, p_{k,2\ell}), \qquad p_{k,2\ell} \equiv \frac{\int_{I_{k,2\ell}} f}{\int_{I_{k-1,\ell}} f} = \frac{f_{k,2\ell}}{f_{k,2\ell} + f_{k,2\ell+1}},$$

so that $W_{k,2\ell}$ is distributed exactly according to $N(0, \sigma_{k-1}^2)$ for $p_{k,2\ell} = \frac{1}{2}$, compared with (2.4). Thus, the distributions of $W_{k,2\ell}$ and $W_{k,2\ell}^\star$ are close at high resolution levels as long as $f$ is sufficiently smooth, even for small $N_{k-1,\ell} = m$.

The equivalence mappings $T_n$, with randomization through $\widetilde{\mathbf{U}}$, are defined by

$$T_n : \{N, \mathbf{X}_N, \widetilde{\mathbf{U}}\} \to \mathbf{W}_{[k_0, \infty)} \to Z_n \equiv \{Z_n(t) : 0 \leq t \leq 1\},$$

where for $k_0 \leq k \leq \infty$, $\mathbf{W}_{[k_0,k)} \equiv \{\overline{\mathbf{Z}}_{k_0}, \mathbf{W}_j, k_0 < j < k\}$, and $\overline{\mathbf{Z}}_k$ and $\mathbf{W}_k$ are as in (2.6) and (2.9). The inverse of $T_n$ is a deterministic many-to-one mapping defined by

$$T_n^{-1} : Z_n^\star \to \mathbf{W}_{[k_0,\infty)}^\star \to (N^\star, \mathbf{X}_{N^\star}^\star),$$

where for $k_0 \leq k \leq \infty$, $\mathbf{W}_{[k_0,k)}^\star \equiv \{\overline{\mathbf{Z}}_{k_0}^\star, \mathbf{W}_j^\star, k_0 < j < k\}$.

REMARK 1. One need only carry out the above construction to $k = k_1 : 2^{k_1} > \varepsilon n$ since we shall assume that $f \in B(\frac{1}{2}, 2, 2)$ and then the observations $\mathbf{W}_{[k_0,k)}^\star \equiv \{\overline{\mathbf{Z}}_{k_0}^\star, \mathbf{W}_j^\star, k_0 < j < k\}$ and $\mathbf{W}_{[k_0,k)} \equiv \{\overline{\mathbf{Z}}_{k_0}, \mathbf{W}_j, k_0 < j < k\}$ are asymptotically sufficient for the Gaussian process and Poisson process experiments. See Brown and Low (1996) for a detailed argument in the context of nonparametric regression.

*Mappings for the density estimation model.* The constructive asymptotic equivalence between density estimation experiments and Gaussian experiments is established by first randomizing the density estimation experiment to an approximation of the Poisson process and then applying the randomized mapping as given above. Set $\gamma_k = \sup_{f \in \mathcal{H}} \|f - \bar{f}_k\|_{1/2,2,2}^2$ and note that since $\mathcal{H}$ is compact in $B(1/2, 2, 2)$, $\gamma_k \downarrow 0$. Now let $k_0$ be the smallest integer such that $4^{k_0}/n \geq \gamma_{k_0}$ and divide the unit interval into subintervals of equal length with length equal to $2^{-k_0}$. Let $\tilde{f}_n$ be the corresponding histogram estimate based on $\mathbf{V}_n^\star$. Now note that since functions $f \in \mathcal{H}$ are bounded below by $\varepsilon_0 > 0$ it follows that

$$(2.11) \quad \int_0^1 (\sqrt{\tilde{f}_n} - \sqrt{f})^2 \leq \int_0^1 (\sqrt{\tilde{f}_n} - \sqrt{f})^2 \frac{(\sqrt{\tilde{f}_n} + \sqrt{f})^2}{\varepsilon_0} = \int_0^1 \frac{(\tilde{f}_n - f)^2}{\varepsilon_0}.$$

Now

$$(2.12) \quad E \int_0^1 (\tilde{f}_n - f)^2 = E \int_0^1 (\tilde{f}_n - \bar{f}_{k_0})^2 + \int_0^1 (f - \bar{f}_{k_0})^2$$

and simple calculations show that the histogram estimate $\tilde{f}_n$ satisfies $E\tilde{f}_n(x) = \bar{f}_{k_0}(x)$ and $\operatorname{Var} \tilde{f}_n(x) \leq \bar{f}_{k_0}(x)\frac{2^{k_0}}{n}$. Hence,

$$(2.13) \quad n^{1/2} E \int_0^1 (\tilde{f}_n - \bar{f}_{k_0})^2 \leq n^{1/2} \frac{2^{k_0}}{n} \leq 2\gamma_{k_0}^{1/2} \to 0.$$

Now $n^{1/2} \leq \frac{2^{k_0}}{\gamma_{k_0}^{1/2}}$ and hence, from (1.7),

$$(2.14) \qquad n^{1/2} \int_0^1 (f - \bar{f}_{k_0})^2 \leq \frac{1}{\gamma_{k_0}^{1/2}} \|f - \bar{f}_{k_0}\|_{1/2,2,2}^2 \leq \gamma_{k_0}^{1/2} \to 0.$$

It thus follows from (2.11) to (2.14) that

$$(2.15) \qquad n^{1/2} \sup_{f \in \mathcal{H}} E \int_0^1 (\sqrt{\tilde{f}_n} - \sqrt{f})^2 \to 0.$$

Hence the density estimate is squared Hellinger consistent at a rate faster than square root of $n$.

Now generate $\tilde{N}$, a Poisson random variable with expectation $n$ and independent of $\mathbf{V}_n^\star$. If $\tilde{N} > n$ generate $\tilde{N} - n$ conditionally independent observations $V_{n+1}^\star, \ldots, V_{\tilde{N}}^\star$ with common density $\tilde{f}_n$. Finally let $(\tilde{N}, \tilde{\mathbf{X}}_{\tilde{N}}) = (\tilde{N}, V_1^\star, V_2^\star, \ldots, V_{\tilde{N}}^\star)$ and write $R_n^1$ for this randomization from $\mathbf{V}_n^\star$ to $(\tilde{N}, \tilde{\mathbf{X}}_{\tilde{N}})$,

$$R_n^1 : \mathbf{V}_n^\star \to (\tilde{N}, \tilde{\mathbf{X}}_{\tilde{N}}).$$

A map from the Poisson number of independent observations back to the fixed number of observations is obtained similarly. This time let $\hat{f}_n$ be the histogram estimator based on $(N, X_N)$. If $N < n$ generate $n - N$ additional conditionally independent observations with common density $\hat{f}_n$. It is also easy to check that

$$(2.16) \qquad n^{1/2} \sup_{f \in \mathcal{H}} E \int_0^1 (\sqrt{\hat{f}_n} - \sqrt{f})^2 \to 0.$$

Now label these observations $\mathbf{V}_n = (V_1, \ldots, V_n)$ and write $R_n^2$ for this randomization from $(N, \mathbf{X}_N)$ to $\mathbf{V}_n$,

$$R_n^2 : (N, \mathbf{X}_N) \to \mathbf{V}_n.$$

REMARK 2. It should also be possible to map the density estimation problem directly into an approximation of the white noise with drift model. Dividing the interval into $2^{k_0}$ subintervals and conditioning on the number of observations falling in each subinterval, the conditional distribution within each subinterval is the same as for the Poisson process. Therefore, it is only necessary to have a version of Theorem 4 for a $2^{k_0}$-dimensional multinomial experiment.

Carter (2002) provides a transformation from a $2^{k_0}$-dimensional multinomial to a multivariate normal as in Theorem 4 such that the total-variation distance between the distributions is $O(k_0 2^{k_0} n^{-1/2})$. The transformation is similar to ours in that it adds uniform noise and then uses the square root as a variance-stabilizing transformation. However, the covariance structure of the multinomial complicates the issue and necessitates using a multi-resolution structure similar to the one applied here to the conditional experiments. The Carter (2002) result can be used in place of Theorem 4 to get a slightly weaker bound on the error in the approximation in Theorem 3 (because of the extra $k_0$ factor) when the total number of observations is fixed. This is enough to establish Theorem 2 if the inequalities bounding $\alpha$ and $\beta$ are changed to strictly greater than. It is also enough to establish Theorem 1 if $\mathcal{H}$ is a Besov space with $\alpha > \frac{1}{2}$. Carter (2000) also showed that a somewhat more complicated transformation leads to a deficiency bound on the normal approximation to the multinomials without the added $k_0$ factor.

**3. Main theorem.** The theorems in Section 1 on the equivalence of white noise with drift experiments and Poisson process experiments are consequences of the following theorem which uniformly bounds the Hellinger distance between the randomized mappings described in Section 2.

THEOREM 3. *Suppose $\inf_{0<x<1} f(x) \geq \varepsilon_0 > 0$. Let $\mathbf{W}^\star_{[k_0,k)} \equiv \{\overline{\mathbf{Z}}^\star_{k_0}, \mathbf{W}^\star_j, k_0 < j < k\}$ with the variables in (2.1) and (2.2), and $\mathbf{W}_{[k_0,k)} \equiv \{\overline{\mathbf{Z}}_{k_0}, \mathbf{W}_j, k_0 < j < k\}$ with the variables in (2.6) and (2.9). Then there exist universal constants $C$, $D_1$ and $D_2$ such that for all $k_1 > k_0$,*

$$H^2(\mathbf{W}^\star_{[k_0,k_1)}, \mathbf{W}_{[k_0,k_1)})$$

$$\leq \frac{C}{\varepsilon_0} \frac{4^{k_0}}{n} + \frac{D_1}{\varepsilon_0^2} \sum_{k=k_0}^{\infty} 2^k \sum_{\ell=0}^{2^k-1} \theta_{k,\ell}^2 + \frac{D_2}{\varepsilon_0^3} \frac{n}{4^{k_0}} \sum_{k=k_0}^{\infty} 2^{3k} \sum_{\ell=0}^{2^k-1} \theta_{k,\ell}^4$$

$$\leq \frac{C}{\varepsilon_0} \frac{4^{k_0}}{n} + \frac{D_1}{\varepsilon_0^2} \|f - \bar{f}_{k_0}\|_{1/2,2,2}^2 + \frac{D_2}{\varepsilon_0^3} \frac{n}{4^{k_0}} \|f - \bar{f}_{k_0}\|_{1/2,4,4}^4,$$

*where $\theta_{k,\ell}$ are the Haar coefficients of $f$ as in (1.3), $\bar{f}_k$ is as in (1.6) and $\|\cdot\|_{1/2,p,p}$ are the Besov norms in (1.5).*

REMARK 3. Here the universal constant $C$ is the same as the one in Theorem 4, while $D_1 = \frac{3D}{8} + 2$ and $D_2 = \frac{D}{9} + \frac{8}{3}$ for the $D$ in Theorem 5.

The proof of Theorem 3 is based on the inequalities established in Sections 4 and 5 for the normal approximation of Poisson and Binomial variables. Some additional technical lemmas are given in the Appendix.

Let $\widetilde{X}_{m,p}$ be a $\text{Bin}(m,p)$ variable, $\widetilde{X}_\lambda$ be a Poisson variable with mean $\lambda$, and $\widetilde{U}$ be a uniform variable in $[-\frac{1}{2}, \frac{1}{2})$ independent of $\widetilde{X}_{m,p}$ and $\widetilde{X}_\lambda$. Define

$$(3.1) \qquad \tilde{g}_{m,p}(x) \equiv \frac{d}{dx} P\{\Phi^{-1}(F_m(\widetilde{X}_{m,p} + \widetilde{U})) \leq x\}$$

with the $F_m$ in (2.8) and the $N(0,1)$ distribution function $\Phi$, and define

(3.2) $$\tilde{g}_\lambda(x) \equiv \frac{d}{dx} P\{2\,\mathrm{sgn}(\widetilde{X}_\lambda + \widetilde{U})\sqrt{|\widetilde{X}_\lambda + \widetilde{U}|} \leq x\}.$$

Write $\varphi_b$ for the density of $N(b,1)$ variables.

PROOF OF THEOREM 3. Let $g^\star_{[k_0,k)}(\mathbf{w}_{[k_0,k)})$ and $g_{[k_0,k)}(\mathbf{w}_{[k_0,k)})$ be the joint densities of $\mathbf{W}^\star_{[k_0,k)}$ and $\mathbf{W}_{[k_0,k)}$, $g^\star_k(\mathbf{w}_k)$ be the joint density of $\mathbf{W}^\star_k$, and $g_k(\mathbf{w}_k|\mathbf{w}_{[k_0,k)})$ be the conditional joint density of $\mathbf{W}_k$ given $\mathbf{W}_{[k_0,k)}$. Since $\mathbf{W}^\star_k$ is independent of $\mathbf{W}^\star_{[k_0,k)}$,

$$\sqrt{g^\star_{[k_0,k)}g_{[k_0,k)}} - \sqrt{g^\star_{[k_0,k+1)}g_{[k_0,k+1)}} = \sqrt{g^\star_{[k_0,k)}g_{[k_0,k)}}(1 - \sqrt{g^\star_k g_k}),$$

so that the Hellinger distance can be written as

(3.3)
$$\begin{aligned} H^2_f(\mathbf{W}^\star_{[k_0,k_1)}, \mathbf{W}_{[k_0,k_1)}) \\ &= 2\left(1 - \int \sqrt{g^\star_{[k_0,k_1)}g_{[k_0,k_1)}}\right) \\ &= 2\left(1 - \int \sqrt{g^\star_{[k_0,k_0+1)}g_{[k_0,k_0+1)}}\right) \\ &\quad + \sum_{k_0 < k < k_1} 2\int \sqrt{g^\star_{[k_0,k)}g_{[k_0,k)}}\left(1 - \int \sqrt{g^\star_k g_k}\right) \\ &= H^2_f(\overline{\mathbf{Z}}^\star_{k_0}, \overline{\mathbf{Z}}_{k_0}) + \sum_{k_0 < k < k_1} \int \sqrt{g^\star_{[k_0,k)}g_{[k_0,k)}} H^2(g^\star_k, g_k). \end{aligned}$$

At the initial resolution level $k_0$, $N_{k_0,\ell}$ are independent Poisson variables by (2.5), so that $\overline{Z}_{k_0,\ell}$ are independent. This and the independence of $\overline{Z}^\star_{k_0,\ell}$ from (2.1) imply

$$H^2_f(\overline{\mathbf{Z}}^\star_{k_0}, \overline{\mathbf{Z}}_{k_0}) \leq \sum_{\ell=0}^{2^{k_0}-1} H^2_f(\overline{Z}^\star_{k_0,\ell}, \overline{Z}_{k_0,\ell}).$$

By (2.6) and (3.2) $\overline{Z}_{k_0,\ell}/\sigma_{k_0}$ have densities $\tilde{g}_{\lambda_{k_0,\ell}}$, while $\overline{Z}^\star_{k_0,\ell}/\sigma_{k_0}$ are $N(h_{k_0,\ell}/\sigma_{k_0}, 1)$ variables by (2.3). Thus, Theorem 4 can be used to obtain

$$H^2_f(\overline{Z}^\star_{k_0,\ell}, \overline{Z}_{k_0,\ell}) = H^2_f(\tilde{g}_{\lambda_{k_0,\ell}}, \varphi_{h_{k_0,\ell}/\sigma_{k_0}}) \leq \frac{C}{\lambda_{k_0,\ell}} + \frac{1}{2}\left(2\sqrt{\lambda_{k_0,\ell}} - \frac{h_{k_0,\ell}}{\sigma_{k_0}}\right)^2.$$

Since $\lambda_{k,\ell} = f_{k,\ell}/(4\sigma_k^2)$ by (2.7) and $\sigma_k^2 = 2^{k-2}/n$ by (2.3), the above calculation yields

$$H^2_f(\overline{\mathbf{Z}}^\star_{k_0}, \overline{\mathbf{Z}}_{k_0}) \leq C \sum_{\ell=0}^{2^{k_0}-1} \frac{2^{k_0}}{nf_{k_0,\ell}} + \sum_{\ell=0}^{2^{k_0}-1} \frac{2n}{2^{k_0}}(\sqrt{f_{k_0,\ell}} - h_{k_0,\ell})^2$$

(3.4)
$$\leq C\frac{2^{2k_0}}{n\varepsilon_0} + \sum_{\ell=0}^{2^{k_0}-1} \frac{n2^{k_0}}{2\varepsilon_0^3}\left(\int_{I_{k_0,\ell}}(f-f_{k_0,\ell})^2\right)^2$$

by Lemma 1(i) and the bound $f \geq \varepsilon_0$.

For $k > k_0$ and $0 \leq \ell < 2^{k-1} - 1$, define

(3.5) $\quad \mu_{k,2\ell} \equiv \sqrt{m_{k,2\ell}}(2p_{k,2\ell}-1), \quad \beta_{k,2\ell} \equiv \sqrt{\lambda_{k-1,\ell}}(2p_{k,2\ell}-1),$

where $p_{k,2\ell}$ are as in (2.10), $\lambda_{k,\ell} = f_{k,\ell}n/2^k$ are as in (2.7), and the functions $m_{k,2\ell} \equiv m_{k,2\ell}(\mathbf{w}_{[k_0,k)})$ are defined by $N_{k-1,\ell} = m_{k,2\ell}(\mathbf{W}_{[k_0,k)})$. At a fixed resolution level $k > k_0$, and for $\ell = 0, \ldots, 2^{k-1}-1$, $N_{k,2\ell}$ are independent binomial variables conditionally on $\mathbf{W}_{[k_0,k)}$, so that by (2.9) and (3.1) $W_{k,2\ell}/\sigma_{k-1}$ are independent variables with densities $\tilde{g}_{m_{k,2\ell},p_{k,2\ell}}$ under the conditional density $g_k$. In addition, $W^*_{k,2\ell}$ are independent normal variables with variance $\sigma^2_{k-1}$ under $g^*_k$. Thus,

(3.6) $$H^2(g^*_k, g_k) \leq \sum_{\ell=0}^{2^{k-1}-1} H^2(\tilde{g}_{m_{k,2\ell},p_{k,2\ell}}, \varphi_{\beta^\star_{k,2\ell}}),$$

by (2.4), where $\beta^\star_{k,2\ell} \equiv EW^\star_{k,2\ell}/\sigma_{k-1} = \sqrt{4n}\int h\phi_{k-1,\ell}$. It follows from Theorem 5 and (3.5) that for fixed $\mathbf{w}_{[k_0,k)}$,

(3.7)
$$H^2(\tilde{g}_{m_{k,2\ell},p_{k,2\ell}}, \varphi_{\beta^\star_{k,2\ell}})$$
$$\leq D\left\{\left[p_{k,2\ell}-\frac{1}{2}\right]^2 + m_{k,2\ell}\left[p_{k,2\ell}-\frac{1}{2}\right]^4\right\} + \frac{(\mu_{k,2\ell}-\beta^\star_{k,2\ell})^2}{2}.$$

Furthermore, it follows from Lemma 3 that

$$\int \sqrt{g^\star_{[k_0,k)}g_{[k_0,k)}}(\sqrt{m_{k,2\ell}} - \sqrt{\lambda_{k-1,\ell}})^2$$
$$\leq \sqrt{\int g_{[k_0,k)}(\sqrt{m_{k,2\ell}} - \sqrt{\lambda_{k-1,\ell}})^4}$$
$$= \sqrt{E(\sqrt{N_{k-1,\ell}} - \sqrt{\lambda_{k-1,\ell}})^4} \leq 2,$$

so that by (3.5),

$$\int \sqrt{g^\star_{[k_0,k)}g_{[k_0,k)}}(\mu_{k,2\ell}-\beta^\star_{k,2\ell})^2 \leq 4(2p_{k,2\ell}-1)^2 + 2(\beta_{k,2\ell}-\beta^\star_{k,2\ell})^2.$$

Similarly, $\int \sqrt{g^\star_{[k_0,k)}g_{[k_0,k)}}m_{k,2\ell} \leq \sqrt{EN^2_{k-1,\ell}} \leq \lambda_{k-1,\ell} + 1/2$. Thus, by (3.7),

$$\int \sqrt{g^\star_{[k_0,k)}g_{[k_0,k)}}H^2(\tilde{g}_{m_{k,2\ell},p_{k,2\ell}}, \varphi_{\beta^\star_{k,2\ell}})$$

(3.8)
$$\leq 4D_1[p_{k,2\ell}-\tfrac{1}{2}]^2 + D\lambda_{k-1,\ell}[p_{k,2\ell}-\tfrac{1}{2}]^4 + (\beta_{k,2\ell}-\beta^\star_{k,2\ell})^2,$$

with $D_1 = 3D/8 + 2$. Now, by (2.10) and (1.3),

$$(3.9) \quad p_{k,2\ell} - \frac{1}{2} = \frac{\int_{I_{k,2\ell}} f - \int_{I_{k,2\ell+1}} f}{2\int_{I_{k-1,\ell}} f} = \frac{\sqrt{2^{k-1}}\theta_{k-1,\ell}}{2f_{k-1,\ell}},$$

so that by (3.5), (2.7), the definition of $\beta^\star_{k,2\ell}$ in (3.6) and Lemma 1(ii),

$$(3.10) \quad \begin{aligned} |\beta_{k,2\ell} - \beta^\star_{k,2\ell}| &= \left|\sqrt{\frac{nf_{k-1,\ell}}{2^{k-1}}} \frac{\sqrt{2^{k-1}}\theta_{k-1,\ell}}{f_{k-1,\ell}} - \sqrt{4n}\int h\phi_{k-1,\ell}\right| \\ &= \sqrt{4n}\left|\frac{\theta_{k-1,\ell}}{2\sqrt{f_{k-1,\ell}}} - \int h\phi_{k-1,\ell}\right| \\ &\leq \sqrt{4n}2^{(k-1)/2-1}f_{k-1,\ell}^{-3/2}\int_{I_{k-1,\ell}}(f-f_{k-1,\ell})^2. \end{aligned}$$

Inserting (3.9) and (3.10) into (3.8) and summing over $\ell$ via (3.6), we find

$$\begin{aligned} & \int \sqrt{g^\star_{[k_0,k)}g_{[k_0,k)}}H^2(g^*_k, g_k) \\ &\leq \sum_{\ell=0}^{2^{k-1}-1} \int \sqrt{g^\star_{[k_0,k)}g_{[k_0,k)}}H^2(\tilde{g}_{m_{k,2\ell},p_{k,2\ell}}, \varphi_{\beta^\star_{k,2\ell}}) \\ (3.11) \quad &\leq \sum_{\ell=0}^{2^{k-1}-1}\left[4D_1\frac{2^k\theta^2_{k-1,\ell}}{8f^2_{k-1,\ell}} + D\lambda_{k-1,\ell}\frac{4^k\theta^4_{k-1,\ell}}{64f^4_{k-1,\ell}}\right. \\ &\qquad\left. + \frac{n2^k}{2f^3_{k-1,\ell}}\left(\int_{I_{k-1,\ell}}(f-f_{k-1,\ell})^2\right)^2\right] \\ &\leq \sum_{\ell=0}^{2^{k-1}-1}\left[\frac{D_1}{\varepsilon_0^2}2^{k-1}\theta^2_{k-1,\ell} + \left(\frac{D}{16}+1\right)\frac{n2^{k-1}}{\varepsilon_0^3}\left(\int_{I_{k-1,\ell}}(f-f_{k-1,\ell})^2\right)^2\right], \end{aligned}$$

due to $\lambda_{k,\ell} = nf_{k,\ell}/2^k$ in (2.7) and $\theta^2_{k,\ell} \leq \int_{I_{k,\ell}}(f-f_{k,\ell})^2$.

Finally, inserting (3.4) and (3.11) into (3.3) and then using Lemma 2 yields

$$\begin{aligned} &H^2_f(\mathbf{W}^\star_{[k_0,k_1+1)}, \mathbf{W}_{[k_0,k_1+1)}) \\ &\leq C\frac{2^{2k_0}}{n\varepsilon_0} + \frac{D_1}{\varepsilon_0^2}\sum_{k=k_0}^{k_1-2}2^k\sum_{\ell=0}^{2^k-1}\theta^2_{k,\ell} \\ &\quad + \left(\frac{D}{16}+\frac{3}{2}\right)\sum_{k=k_0}^{k_1-2}\sum_{\ell=0}^{2^k-1}\frac{n2^k}{\varepsilon_0^3}\left(\int_{I_{k,\ell}}(f-f_{k,\ell})^2\right)^2 \end{aligned}$$

$$\le \frac{C}{\varepsilon_0}\frac{4^{k_0}}{n} + \frac{D_1}{\varepsilon_0^2}\sum_{k=k_0}^{k_1-2} 2^k \sum_{\ell=0}^{2^k-1}\theta_{k,\ell}^2 + \frac{D_2}{\varepsilon_0^3}\frac{n}{4^{k_0}}\sum_{k=k_0}^{\infty} 2^{3k}\sum_{\ell=0}^{2^k-1}\theta_{k,\ell}^4,$$

with $D_2 \equiv (\frac{D}{16} + \frac{3}{2})/(1-\frac{1}{4})^2 = \frac{D}{9} + \frac{8}{3}$ and the theorem follows. □

**4. Approximation of Poisson variables.** Let $\widetilde{X}_\lambda$ be a Poisson random variable with mean $\lambda$ and $\widetilde{U}$ be a uniform variable on $[-\frac{1}{2},\frac{1}{2})$ independent of $\widetilde{X}_\lambda$. Define

(4.1) $\quad \widetilde{Z}_\lambda \equiv 2\,\mathrm{sgn}(\widetilde{X}_\lambda + \widetilde{U})\sqrt{|\widetilde{X}_\lambda + \widetilde{U}|}, \qquad \tilde{g}_\lambda(y) \equiv \frac{d}{dy}P\{\widetilde{Z}_\lambda \le y\}.$

The main result of this section is a local limit theorem which bounds the squared Hellinger distance between this transformed Poisson random variable and a normal random variable.

THEOREM 4. *Let $\widetilde{Z}_\lambda$ and $\tilde{g}_\lambda$ be as in (4.1). Let $Z_\lambda^* \sim N(2\sqrt{\lambda},1)$ and $\varphi_\mu$ be the density of $N(0,\mu)$. Let $H(\cdot,\cdot)$ be the Hellinger distance. Then, as $\lambda \to \infty$,*

(4.2) $\qquad H^2(\widetilde{Z}_\lambda, Z_\lambda^*) = H^2(\tilde{g}_\lambda, \varphi_{2\sqrt{\lambda}}) = (7+o(1))\frac{1}{96\lambda}.$

*Consequently, there exists a universal constant $C < \infty$ such that*

(4.3) $\qquad H^2(\tilde{g}_\lambda, \varphi_\mu) \le C/\lambda + (2\sqrt{\lambda} - \mu)^2/2 \qquad \forall \lambda > 0, \mu.$

REMARK 4. The theorem remains valid if $\widetilde{Z}_\lambda$ is replaced by

$$\widetilde{Z}'_\lambda \equiv 2\sqrt{\widetilde{X}_\lambda + \widetilde{U} + \tfrac{1}{2}},$$

since $H^2(\widetilde{Z}_\lambda, \widetilde{Z}'_\lambda)$ is bounded by

$$2 - 2\int \sqrt{f_{|\widetilde{X}_\lambda + \widetilde{U}|}f_{\widetilde{X}_\lambda + \widetilde{U} + 1/2}} \le 2 - \left\{1 + \sum_{j=0}^{\infty} e^{-\lambda}\sqrt{\frac{\lambda^{2j+1}}{j!(j+1)!}}\right\}$$

$$= 1 - E\sqrt{\frac{\widetilde{X}_\lambda}{\lambda}} \le \min\left(1, \frac{C'}{\lambda}\right).$$

PROOF OF THEOREM 4. The second inequality of (4.3) follows immediately from (4.2), since $H^2(\varphi_{\mu_1}, \varphi_{\mu_2}) = (\mu_1 - \mu_2)^2/4$ [cf. Brown, Cai, Low and Zhang (2002), Lemma 3] and $H^2(\tilde{g}_\lambda, \varphi_\mu) \le 2$.

Let $t(x) \equiv 2\,\mathrm{sgn}(x)\sqrt{|x|}$, a strictly increasing function. Define

(4.4) $\qquad \widetilde{X}_\lambda^* \equiv t^{-1}(Z_\lambda^*) = \mathrm{sgn}(Z_\lambda^*)(Z_\lambda^*)^2/4.$

Let $f_\lambda$ and $f_\lambda^*$ denote the densities of $\widetilde{X}_\lambda + \widetilde{U}$ and $\widetilde{X}^*$, respectively. Since $t(\cdot)$ is invertible, $H^2(\widetilde{Z}_\lambda, Z_\lambda^*) = H(\widetilde{X}_\lambda + \widetilde{U}, \widetilde{X}_\lambda^*) = 2 - 2\int \sqrt{f_\lambda f_\lambda^*}$, so that it suffices to show

$$(4.5) \qquad A_\lambda \equiv \int \sqrt{f_\lambda f_\lambda^*} = 1 - \frac{C_\lambda}{\lambda}, \qquad \lim_{\lambda \to \infty} C_\lambda = \tfrac{7}{192}.$$

Since $\widetilde{U}$ is uniform, $f_\lambda(x) = e^{-\lambda}\lambda^j/j!$ on $[j - 1/2, j + 1/2)$, so that

$$(4.6) \qquad A_\lambda = \sum_{j=0}^{\infty} f_\lambda(j) \int_{j-1/2}^{j+1/2} \{f_\lambda^*(x)/f_\lambda(j)\}^{1/2} dx.$$

Since $t'(x) = |x|^{-1/2}$, by (4.4) $f_\lambda^*(x) = |x|^{-1/2} \varphi(t(x) - 2\sqrt{\lambda})$. This gives

$$\frac{f_\lambda^*(x)}{f_\lambda(j)} = \frac{\exp\{-(2\sqrt{x} - 2\sqrt{\lambda})^2/2\}}{\sqrt{2\pi x} e^{-\lambda}\lambda^j/j!} = \exp[2\psi_j(x)], \qquad j - \tfrac{1}{2} \leq x < j + \tfrac{1}{2},$$

for $j \geq 1$, in view of the Stirling formula $j! = \sqrt{2\pi} j^{j+1/2} \exp(-j + \varepsilon_j)$, where

$$(4.7) \quad \psi_j(x) \equiv -(\sqrt{x} - \sqrt{\lambda})^2 - \frac{\log x}{4} + \frac{\lambda}{2} + \frac{j}{2}\log\left[\frac{j}{\lambda}\right] + \frac{\log j}{4} - \frac{j}{2} + \frac{\varepsilon_j}{2}$$

with $1/(12j + 1) < \varepsilon_j < 1/(12j)$, for $j = 1, 2, \ldots$. Now, by the mean-value theorem,

$$\int_{j-1/2}^{j+1/2} \left\{\frac{f_\lambda^*(x)}{f_\lambda(j)}\right\}^{1/2} dx$$
$$= \int_{j-1/2}^{j+1/2} \exp\left[\psi_j(j) + \psi_j'(j)(x - j) + \frac{\psi_j''(x_j)}{2}(x - j)^2\right] dx$$

for some $|x_j - j| \leq \tfrac{1}{2}$, with

$$(4.8) \qquad \psi_j'(x) = \sqrt{\frac{\lambda}{x}} - 1 - \frac{1}{4x}, \qquad \psi_j''(x) = -\frac{\sqrt{\lambda}}{2x^{3/2}} + \frac{1}{4x^2}.$$

Since $\exp[\psi_j(j) + \psi_j''(x_j)(x - j)^2/2]$ is symmetric about $j$, it follows that

$$(4.9) \quad \int_{j-1/2}^{j+1/2} \left\{\frac{f_\lambda^*(x)}{f_\lambda(j)}\right\}^{1/2} dx$$
$$= \int_{j-1/2}^{j+1/2} \exp\left[\psi_j(j) + \frac{\psi_j''(x_j)(x-j)^2}{2}\right] \sum_{k=0}^{\infty} \frac{(\psi_j'(j)(x-j))^{2k}}{(2k)!} dx.$$

Now, we shall take uniform Taylor expansions of $\psi_j$ and their derivatives in

$$J_\lambda \equiv \{j : |j/\lambda - 1| \leq \lambda^{-2/5}\}.$$

By (4.7), $\psi_j(j) = \lambda \psi(j/\lambda) + \varepsilon_j/2$ with

$$\psi(x) \equiv -(\sqrt{x} - 1)^2 + \frac{1-x}{2} + \frac{x}{2} \log x.$$

Since $\psi(1) = \psi'(1) = \psi''(1) = 0$, $\psi'''(1) = 1/4$ and $\psi''''(1) = -\frac{7}{8}$,

$$\lambda \psi\left(\frac{j}{\lambda}\right) = \frac{\lambda}{4} \frac{(j-\lambda)^3}{3!\lambda^3} - \frac{7\lambda}{8} \frac{(j-\lambda)^4}{4!\lambda^4} (1 + o(1)) = o(1).$$

Since $1/(12j+1) < \varepsilon_j < 1/(12j)$, $\varepsilon_j/2 = (1+o(1))/(24\lambda) = o(1)$. Thus,

$$\psi_j(j) = \frac{(j-\lambda)^3}{24\lambda^2} - \frac{7}{8} \frac{(j-\lambda)^4}{24\lambda^3} (1 + o(1)) + \frac{1 + o(1)}{24\lambda} = o(1)$$

uniformly in $J_\lambda$ as $\lambda \to \infty$. Similarly, by (4.8) and $|x_j - j| \leq \frac{1}{2}$,

$$\{\psi'_j(j)\}^2 = (1 + o(1)) \frac{(j-\lambda)^2}{4\lambda^2} + \frac{o(1)}{\lambda} = o(1),$$

$$\psi''_j(x_j) = \frac{-1 + o(1)}{2\lambda} = o(1).$$

These expansions and (4.9) imply that uniformly in $J_\lambda$,

$$\int_{j-1/2}^{j+1/2} \left\{\frac{f^*_\lambda(x)}{f_\lambda(j)}\right\}^{1/2} dx$$

$$= \int_{j-1/2}^{j+1/2} \left[1 + \psi_j(j) + \{\psi''_j(x_j) + (\psi'_j(j))^2\} \frac{(x-j)^2}{2}\right] dx$$

$$+ o(1) \sum_{k=0}^{2} \frac{(j-\lambda)^{2k}}{\lambda^{k+1}}$$

$$= 1 + \frac{(j-\lambda)^3}{24\lambda^2} - \frac{7}{8} \frac{(j-\lambda)^4}{24\lambda^3} + \frac{1}{24\lambda} + \frac{1}{24}\left[\frac{-1}{2\lambda} + \frac{(j-\lambda)^2}{4\lambda^2}\right]$$

$$+ o(1) \sum_{k=0}^{2} \frac{(j-\lambda)^{2k}}{\lambda^{k+1}},$$

as $\int_{j-1/2}^{j+1/2} (x-j)^2 dx = \frac{1}{12}$. Since $f_\lambda(j)$ is the Poisson probability mass function of $\widetilde{X}_\lambda$,

$$\sum_{j \in J_\lambda} f_\lambda(j) \int_{j-1/2}^{j+1/2} \left\{\frac{f^*_\lambda(x)}{f_\lambda(j)}\right\}^{1/2} dx$$

(4.10)
$$= 1 + \frac{1}{24\lambda} - \left[\frac{7}{8}\right] \frac{3}{24\lambda} + \frac{1}{24\lambda} - \frac{1}{96\lambda} + \frac{o(1)}{\lambda} = 1 - \frac{7 + o(1)}{192\lambda}$$

as $\sum_{j \in J_\lambda} f_\lambda(j) = 1 + o(1/\lambda)$. Note that $E(\widetilde{X}_\lambda - \lambda)^3 = \lambda$ and $E(\widetilde{X}_\lambda - \lambda)^4 = 3\lambda^2 + \lambda$. Hence, (4.5) follows from (4.6), (4.10) and the fact that

$$\sum_{j \notin J_\lambda} f_\lambda(j) \int_{j-1/2}^{j+1/2} \left\{ \frac{f_\lambda^*(x)}{f_\lambda(j)} \right\}^{1/2} dx \leq \sqrt{P\{\widetilde{X}_\lambda \notin J_\lambda\} P\{\widetilde{X}_\lambda^* \notin J_\lambda\}} = o\!\left(\frac{1}{\lambda}\right). \quad \square$$

**5. Approximation of binomial variables.** The strong approximation of a normal by a binomial depends on the cumulative distribution function $F_m$ in (2.8). The addition of the independent uniform $\widetilde{U}$ in (2.8) to the binomial $\widetilde{X}_{m,1/2}$ makes the c.d.f. continuous and thus $\Phi^{-1} \circ F_m$ is a one-to-one function on $(-\frac{1}{2}, m + \frac{1}{2})$ that maps symmetric binomials to standard normals.

Let $\varphi_b$ be the $N(b, 1)$ density and $\tilde{g}_{m,p}$ be the probability density of

(5.1) $$\Phi^{-1}(F_m[\widetilde{X}_{m,p} + \widetilde{U}]), \qquad \widetilde{X}_{m,p} \sim \text{Bin}(m, p),$$

as in (3.1), where $\widetilde{U}$ is an independent uniform on $[-\frac{1}{2}, \frac{1}{2})$.

THEOREM 5.    *There is a constant $C_1 > 0$ such that, for all $m \geq 0$,*

(5.2) $$H^2(\tilde{g}_{m,p}, \varphi_b) = \int (\sqrt{\tilde{g}_{m,p}} - \sqrt{\varphi_b})^2 \, dz \leq C_1 \left( \frac{b^2}{m} + \frac{b^8}{m^2} \right),$$

*where $b = (\sqrt{m}/2) \log(p/(1-p))$. Consequently,*

(5.3) $$H^2(\tilde{g}_{m,p}, \varphi_\beta) \leq D\left[ \left(p - \frac{1}{2}\right)^2 + m\left(p - \frac{1}{2}\right)^4 \right] + \frac{(\sqrt{m}(2p-1) - \beta)^2}{2}.$$

PROOF.    The case when $m = 0$ is trivial because $X = 0$ with probability 1 and therefore $\tilde{g}_{0,p}$ is exactly an $\mathcal{N}(0, 1)$. Thus, the following assumes that $m \geq 1$.

It follows from (3.1) that

(5.4) $$\tilde{g}_{m,p}(z) = p^j (1-p)^{m-j} 2^m \varphi_0(z),$$

where $j = j(z)$ is the integer between 0 and $m$ such that

(5.5) $$\Phi^{-1}[F_m(j - \tfrac{1}{2})] \leq z < \Phi^{-1}[F_m(j + \tfrac{1}{2})].$$

Let $\theta = \log(p/q)$ so that

$$\log \frac{g_{m,p}(z)}{\varphi_0(z)} = \theta\left(j - \frac{m}{2}\right) + \frac{m \log(4pq)}{2},$$

and the second term can be approximated by

(5.6) $$-\frac{\theta^2}{4} - \frac{\theta^4}{24} \leq \log(4pq) = -\log\left[\frac{2 + e^\theta + e^{-\theta}}{4}\right] \leq -\frac{\theta^2}{4} + \frac{\theta^4}{32}.$$

Let $h_1(\theta) = (2 + e^{-\theta} + e^{-\theta})/4$. The second inequality in (5.6) follows from $\log(h_1(\theta)) \geq \log(1 + \theta^2/4) \geq \theta^2/4 - \theta^4/32$. The first inequality in (5.6) follows from $h_1(\theta) \leq 1 + \theta^2/4 + \theta^4/24$ for $|\theta| \leq 4$, and from $\log(h_1(\theta)) \leq |\theta| \leq \theta^2/4$ for $|\theta| > 4$. Now, let

$$(5.7) \qquad z' = z'(z) = \frac{j(z) - m/2}{\sqrt{m}/2} \quad \text{and} \quad b = \theta\frac{\sqrt{m}}{2}.$$

Then for some $-1/24 \leq h_2(\theta) \leq 1/32$ the log ratio is

$$\log \frac{\tilde{g}_{m,p}(z)}{\varphi_0(z)} = z'b - \frac{b^2}{2} + h_2(\theta)m\theta^4.$$

The log ratio of normals with different means is $\log(\varphi_0/\varphi_b) = -zb + b^2/2$. Therefore the ratio with respect to the normal with mean $b$ is

$$(5.8) \qquad \log \frac{\tilde{g}_{m,p}}{\varphi_b} = h_2(\theta)m\theta^4 - b(z - z'), \qquad |h_2(\theta)| \leq \tfrac{1}{24}.$$

Since $y\log(x/y) \leq x - y \leq x\log(x/y)$, for all positive $x$ and $y$,

$$\frac{1}{2}\sqrt{\tilde{g}_{m,p}}\log\left(\frac{\varphi_b}{\tilde{g}_{m,p}}\right) \leq \sqrt{\varphi_b} - \sqrt{\tilde{g}_{m,p}} \leq \frac{1}{2}\sqrt{\varphi_b}\log\left(\frac{\varphi_b}{\tilde{g}_{m,p}}\right),$$

so that by (5.8),

$$(5.9) \qquad \begin{aligned} H^2(\tilde{g}_{p,m}, \varphi_b) &\leq \frac{1}{4}\int\left\{\log\left(\frac{\varphi_b}{\tilde{g}_{m,p}}\right)\right\}^2 (\varphi_b + \tilde{g}_{m,p})\,dz \\ &\leq \left(\frac{m\theta^4}{24}\right)^2 + \frac{b^2}{2}\int (z - z')^2(\varphi_b + \tilde{g}_{m,p})\,dz. \end{aligned}$$

It follows from Carter and Pollard (2004) that the difference between $z$ and $z' = z'(z)$ is bounded by

$$(5.10) \qquad |z - z'| \leq \begin{cases} C_2(m^{-1/2} + m^{-1}|z|^3), & \text{for all } z, \\ C_2(m^{-1/2} + m^{-1}|z'|^3), & \text{if } |z| \leq \sqrt{2m}, \end{cases}$$

for some constant $C_2$. Thus,

$$(5.11) \quad \int (z - z')^2 \tilde{g}_{m,p}\,dz \leq 2C_2^2\left(\frac{1}{m} + \int \frac{|z'|^6}{m^2}\tilde{g}_{m,p}\,dz + \int_{z^2 > 2m} \frac{z^6}{m^2}\tilde{g}_{m,p}\,dz\right).$$

Since $\int \tilde{g}_{m,p} I\{z' = (j - m/2)/\sqrt{m}\}\,dz = P\{\tilde{X}_{m,p} = j\}$,

$$\int |z'|^6 \tilde{g}_{m,p}\,dz = E\left(\frac{\tilde{X}_{m,p} - m/2}{\sqrt{m}}\right)^6 = O\left(1 + m^3\left(\frac{p-1}{2}\right)^6\right) = O(1 + b^6)$$

uniformly in $(m, p)$. It follows from (5.4) that

$$\int_{z^2 > 2m} z^6 \tilde{g}_{m,p}\,dz \leq 2^m \int_{z^2 > 2m} z^6 \varphi_0\,dz = O(2^m m^6 e^{-m}) = O(m^{-1}).$$

The above two inequalities and (5.11) imply

$$\int (z-z')^2 \tilde{g}_{m,p}\, dz \leq 2C_2{}^2 O(1/m + b^6/m^2).$$

Similarly, $\int (z-z')^2 \varphi_b\, dz \leq 2C_2{}^2 O(1/m + b^6/m^2)$. Inserting these two inequalities into (5.9) yields (5.2) in view of (5.7).

Now let us prove (5.3). The Hellinger distance is bounded by 2, so that $b^8/m^2$ in (5.2) can be replaced by $b^4/m$ and it suffices to consider $|p - \frac{1}{2}| \leq \frac{1}{4}$ for the proof of (5.3). By inspecting the infinite series expansion of $\log(\frac{p}{q}) = \log(1+x) - \log(1-x)$ for $x = 2p-1$, we find that for $|p - \frac{1}{2}| \leq \frac{1}{4}$, $|\log(\frac{p}{q})| \leq \frac{8}{3}|2p-1|$ and $|\log(\frac{p}{q}) - 4(p - \frac{1}{2})| \leq \frac{8}{9}|2p-1|^3$. These inequalities, respectively, imply

$$\frac{b^2}{m} + \frac{b^4}{m^2} \leq \frac{16}{9}(2p-1)^2 + \frac{256}{81}m(2p-1)^4$$

and $|b - \sqrt{m}(2p-1)|^2 \leq \frac{16}{81}m|2p-1|^6 \leq \frac{4}{81}m|2p-1|^4$, in view of the definition of $b$, which then imply (5.3) via (5.2) and the fact that $H^2(\varphi_b - \varphi_\beta) = (b-\beta)^2/4$.

□

## APPENDIX

**A.1. The Tusnády inequality.** The coupling of symmetric binomials and normals maps the integers $j$ onto intervals $[\beta_j, \beta_{j+1}]$ such that the normal$(m/2, m/4)$ probability in the interval is equal to the binomial probability at $\binom{m}{j}2^{-j}$. Taking the standardized values

$$z_j = \frac{2(\beta_j - m/2)}{\sqrt{m}}, \qquad u_j = \frac{2(j - 1/2 - m/2)}{\sqrt{m}},$$

Carter and Pollard (2004) showed that for $m/2 < j < m$ and certain universal finite constants $C_\pm$

$$C_- \frac{u_j + 1}{m} \leq z_j - u_j\sqrt{1 + 2\frac{u_j^2}{m}\gamma\left(\frac{u_j}{\sqrt{m}}\right)} - \frac{\log(1 - u_j^2/m)}{2cu_j} \leq C_+ \frac{u_j + \log m}{m}$$

where $c = \sqrt{2\log 2}$ and $\gamma$ is an increasing function with $\gamma(0) = 1/12$ and $\gamma(1) = \log 2 - 1/2$.

This immediately implies that

(A.1) $$|z_j - u_j| \leq \frac{C_0}{m}(|u_j|^3 + \log m) \qquad \forall \frac{u_j^2}{m} \leq \frac{1}{2}$$

for a certain universal constant $C_0 < \infty$. We shall prove (5.10) here based on (A.1). Because of the symmetry in both distributions, it is only necessary to consider $z > 0$.

It follows from (5.5) and (5.7) that

$$z_j \le z < z_{j+1} \iff u_j \le z' = z'(z) < u_{j+1}.$$

Let $z_j \le z < z_{j+1}$. Since $u_{j+1} - u_j = 2/\sqrt{m}$, for $u_{j+1}^2 \le m/2$ (A.1) implies

$$(A.2) \quad |z - z'| \le |z_j - u_j| \vee |z_{j+1} - u_{j+1}| + \frac{2}{\sqrt{m}} \le C_0'\left(\frac{1}{\sqrt{m}} + \frac{|z|^3 \wedge |z'|^3}{m}\right).$$

Since $u_j$ and $z_j$ are both increasing in $j$, it follows that $(z \wedge z')/\sqrt{m}$ are uniformly bounded away from zero for $u_{j+1} \ge \sqrt{m/2}$, so that

$$(A.3) \quad |z - z'| \le |z_j - u_j| \vee |z_{j+1} - u_{j+1}| + \frac{2}{\sqrt{m}} \le C_0'' \frac{|z|^3 \wedge |z'|^3}{m}$$

for $(m+1)/\sqrt{m} = u_{m+1} \ge u_{j+1} \ge m/2$ and $z \le \sqrt{2m}$. Since $z \vee z' \le z \vee u_{m+1} \le \sqrt{2}z$ for $z > \sqrt{2m}$, (A.2) and (A.3) imply

$$|z - z'| \le \begin{cases} C_2(m^{-1/2} + m^{-1}|z|^3), & \text{for all } z, \\ C_2(m^{-1/2} + m^{-1}|z'|^3), & \text{if } |z| \le \sqrt{2m}, \end{cases}$$

for a certain universal $C_2 < \infty$, that is, (5.10).

**A.2. Technical lemmas.** The following three lemmas simplify the rest of the proof of Theorem 3.

LEMMA 1. (i) *Let $f_{k,\ell}$ and $h_{k,\ell}$ be as in (2.7) and (2.3). Then*

$$(A.4) \qquad 0 \le \sqrt{f_{k,\ell}} - h_{k,\ell} \le 2^{k-1} f_{k,\ell}^{-3/2} \int_{I_{k,\ell}} (f - f_{k,\ell})^2.$$

(ii) *Let $\theta_{k,\ell}$ be the Haar coefficients of $f$ as in (1.3). Then*

$$(A.5) \qquad \left|\int h\phi_{k,\ell} - \frac{\theta_{k,\ell}}{2\sqrt{f_{k,\ell}}}\right| \le 2^{k/2-1} f_{k,\ell}^{-3/2} \int_{I_{k,\ell}} (f - f_{k,\ell})^2.$$

PROOF. Let $T = (f - f_{k,\ell})/f_{k,\ell} \ge -1$. By algebra,

$$\sqrt{1+T} - 1 = \frac{T}{1+\sqrt{1+T}} = \frac{T}{2} - \frac{T^2}{2(1+\sqrt{1+T})^2}.$$

It follows from (2.3) and (2.7) that

$$h_{k,\ell} = 2^k \sqrt{f_{k,\ell}} \int_{I_{k,\ell}} \sqrt{1+T}$$

$$= 2^k \sqrt{f_{k,\ell}} \int_{I_{k,\ell}} \left(1 + \frac{f - f_{k,\ell}}{2f_{k,\ell}} - \frac{(f - f_{k,\ell})^2}{2f_{k,\ell}^2(1+\sqrt{1+T})^2}\right),$$

which implies (A.4) as $2^k \int_{I_{k,\ell}} = 1$ and by (2.7) $\int_{I_{k,\ell}} (f - f_{k,\ell}) = 0$. For (ii) we have

$$\int h\phi_{k,\ell} = \sqrt{f_{k,\ell}} \int \phi_{k,\ell} \sqrt{1+T}$$

$$= \sqrt{f_{k,\ell}} \int \phi_{k,\ell} \left(1 + \frac{f - f_{k,\ell}}{2f_{k,\ell}} - \frac{(f - f_{k,\ell})^2}{2f_{k,\ell}^2(1+\sqrt{1+T})^2}\right),$$

which implies (A.5) as $\int \phi_{k,\ell} = 0$ and $|\phi_{k,\ell}| \leq \sqrt{2^k}$ by (1.4). □

LEMMA 2. *Let $\theta_{k,\ell}$ be the Haar coefficients in (1.3) and $f_{k,\ell}$ be as in (2.7). Then*

$$\sum_{k=k_0}^{\infty} 2^k \sum_{\ell=0}^{2^k-1} \left(\int_{I_{k,\ell}} (f - f_{k,\ell})^2\right)^2 \leq \frac{2^{-ck_0}}{(1-1/2^c)^2} \sum_{k=k_0}^{\infty} 2^{k(1+c)} \sum_{\ell=0}^{2^k-1} \theta_{k,\ell}^4 \qquad \forall c > 0.$$

PROOF. Define

$$\delta_{i,j,k,\ell} \equiv \begin{cases} 1, & \text{if } I_{i,j} \subseteq I_{k,\ell}, \\ 0, & \text{otherwise.} \end{cases}$$

Since $\sum_j \delta_{i,j,k,\ell} = 2^{i-k}$ for $i \geq k$, using Cauchy–Schwarz twice yields

$$\left(\int_{I_{k,\ell}} (f - \bar{f}_k)^2\right)^2 = \left(\sum_{i=k}^{\infty} \sum_{j=0}^{2^i-1} \delta_{i,j,k,\ell} \theta_{i,j}^2\right)^2$$

$$\leq \left[\sum_{i=k}^{\infty} 2^{-ic/2} \left(2^{ic} 2^{i-k} \sum_{j=0}^{2^i-1} \delta_{i,j,k,\ell} \theta_{i,j}^4\right)^{1/2}\right]^2$$

$$\leq \sum_{i=k}^{\infty} 2^{-ic} \sum_{i=k}^{\infty} 2^{ic} 2^{i-k} \sum_{j=0}^{2^i-1} \delta_{i,j,k,\ell} \theta_{i,j}^4$$

$$\leq \frac{2^{-k(1+c)}}{1-1/2^c} \sum_{i=k}^{\infty} 2^{i(1+c)} \sum_{j=0}^{2^i-1} \delta_{i,j,k,\ell} \theta_{i,j}^4.$$

Since $\sum_{\ell=0}^{2^k-1} \delta_{i,j,k,\ell} = 1$ for $i \geq k$, the above inequality implies

$$\sum_{k=k_0}^{\infty} 2^k \sum_{\ell=0}^{2^k-1} \left(\int_{I_{k,\ell}} (f - f_{k,\ell})^2\right)^2$$

$$\leq \sum_{k=k_0}^{\infty} 2^k \frac{2^{-k(1+c)}}{1-1/2^c} \sum_{i=k}^{\infty} 2^{i(1+c)} \sum_{j=0}^{2^i-1} \sum_{\ell=0}^{2^k-1} \delta_{i,j,k,\ell} \theta_{i,j}^4$$

$$= \sum_{i=k_0}^{\infty} \left( \sum_{k=k_0}^{i} \frac{2^{-ck}}{1-1/2^c} \right) 2^{i(1+c)} \sum_{j=0}^{2^i-1} \theta_{i,j}^4$$

$$\leq \frac{2^{-ck_0}}{(1-1/2^c)^2} \sum_{i=k_0}^{\infty} 2^{i(1+c)} \sum_{j=0}^{2^i-1} \theta_{i,j}^4. \qquad \square$$

LEMMA 3. *Let $\widetilde{X}_\lambda$ be a Poisson random variable with mean $\lambda$. Then*

$$E(\sqrt{\widetilde{X}_\lambda} - \sqrt{\lambda})^4 \leq 4.$$

PROOF. Since $E(\widetilde{X}_\lambda - \lambda)^4 = \lambda(3\lambda + 1)$,

$$E(\sqrt{\widetilde{X}_\lambda} - \sqrt{\lambda})^4 \leq \frac{E(\widetilde{X}_\lambda - \lambda)^4}{(\sqrt{\lambda}+1)^4} + \lambda^2 P(\widetilde{X}_\lambda = 0)$$

$$\leq \frac{3\lambda+1}{\lambda+6} + 1 \leq 4. \qquad \square$$

**A.3. Proof of Theorem 1.** First note that

$$H(T_n R_n^1 \mathbf{V}_n^\star, \mathbf{Z}_n^\star) \leq H(T_n R_n^1 \mathbf{V}_n^\star, T_n(N, \mathbf{X}_N)) + H(T_n(N, \mathbf{X}_N), \mathbf{Z}_n^\star)$$

and

$$H(\mathbf{V}_n^\star, R_n^2 T_n^{-1} \mathbf{Z}_n^\star) \leq H(\mathbf{V}_n^\star, R_n^2(N, \mathbf{X}_N)) + H(R_n^2(N, \mathbf{X}_N), R_n^2 T_n^{-1} \mathbf{Z}_n^\star).$$

Note also that since for any randomization $T$ and random $X$ and $Y$, $H(TX, TY) \leq H(X, Y)$, it follows that

$$H(T_n R_n^1 \mathbf{V}_n^\star, T_n(N, \mathbf{X}_N)) \leq H(R_n^1 \mathbf{V}_n^\star, (N, \mathbf{X}_N))$$

and

$$H(R_n^2(N, \mathbf{X}_N), R_n^2 T_n^{-1} \mathbf{Z}_n^\star) \leq H((N, \mathbf{X}_N), T_n^{-1} \mathbf{Z}_n^\star) = H(T_n(N, \mathbf{X}_N), \mathbf{Z}_n^\star).$$

For the class $\mathcal{H}$ and the randomizations $R_n^1$ and $R_n^2$ it follows from (2.15), (2.16) and the proof of Proposition 3 on page 508 of Le Cam (1986) that

$$\sup_{f \in \mathcal{H}} H(R_n^1 \mathbf{V}_n^\star, (N, \mathbf{X}_N)) \to 0$$

and

$$\sup_{f \in \mathcal{H}} H(\mathbf{V}_n^\star, R_n^2(N, \mathbf{X}_N)) \to 0.$$

Hence (1.9) and (1.8) will follow once

(A.6) $$\sup_{f \in \mathcal{H}} H(T_n(N, \mathbf{X}_N), \mathbf{Z}_n^\star) \to 0$$

is established.

By Theorem 3, for (A.6) to hold it is sufficient to show that

$$\sup_{f \in \mathcal{H}} \left( \frac{4^{k_0}}{n} + \|f - \bar{f}_{k_0}\|_{1/2,2,2}^2 + \frac{n}{4^{k_0}} \|f - \bar{f}_{k_0}\|_{1/2,4,4}^4 \right) \to 0.$$

If the class of functions $\mathcal{H}$ is a compact set in the Besov spaces, then the partial sums converge uniformly to 0,

$$\sup_{f \in \mathcal{H}} \|f - \bar{f}_k\|_{1/2,p,p} \to 0$$

for $p = 2$ or $4$ as $k \to \infty$. This implies that there is a sequence $\gamma_k \to 0$ such that $\gamma_k^{-1} \sup_{f \in \mathcal{H}} \|f - \bar{f}_k\|_{1/2,4,4}^4 \to 0$. To be specific, let

$$\gamma_k = \sup_{f \in \mathcal{H}} \|f - \bar{f}_k\|_{1/2,4,4}^2.$$

It is necessary to choose the sequence of integers $k_0(n)$ that will be the critical dimension that divides the two techniques. Let $k_0$ be the smallest integer such that $\frac{4^{k_0}}{n} \geq \gamma_{k_0}$. Therefore, $k_0(n) \to \infty$, and as $n \to \infty$,

$$\sup_{f \in \mathcal{H}} \left( \frac{4^{k_0}}{n} + \|f - \bar{f}_{k_0}\|_{1/2,2,2}^2 + \frac{n}{4^{k_0}} \|f - \bar{f}_{k_0}\|_{1/2,4,4}^4 \right)$$

$$\leq \sup_{f \in \mathcal{H}} \left( 4\gamma_{k_0} + \|f - \bar{f}_{k_0}\|_{1/2,2,2}^2 + \frac{1}{\gamma_{k_0}} \|f - \bar{f}_{k_0}\|_{1/2,4,4}^4 \right) \to 0. \quad \square$$

**A.4. Proof of Theorem 2.** Theorem 2 follows from Theorem 1 and the fact that the Lipschitz and Sobolev spaces described are compact in the Besov spaces.

The Lipschitz class is equivalent to $\mathcal{B}_{\beta,\infty,\infty}$ and therefore is compact in $\mathcal{B}_{1/2,p,p}$ if $\beta > \frac{1}{2}$. The Sobolev class is equivalent to $\mathcal{B}_{\alpha,2,2}$ and

$$\|f - \bar{f}_{k_0}\|_{\alpha,2,2}^2 \leq C_\alpha \sum_n |c_n(f)|^2 n^{2\alpha},$$

where $C_\alpha$ depends only on $\alpha$. Thus if $\mathcal{F}$ is compact in Sobolev$(\alpha)$ for $\alpha \geq \frac{1}{2}$ then it is compact in $\mathcal{B}_{1/2,2,2}$.

Further restrictions are required to show that the Sobolev$(\alpha)$ class is compact in $\mathcal{B}_{1/2,4,4}$. If $\|f\|_\beta^{(L)} \leq C_{(L)}$, then $\|\bar{f}_k - \bar{f}_{k+1}\|_\infty \leq C_{(L)} 2^{-k\beta}$, so that

$$\|f - \bar{f}_{k_0}\|_{1/2,4,4}^4 \leq C_{(L)}^2 \sum_{k=k_0}^\infty 2^{k2(1-\beta)} \int |\bar{f}_k - \bar{f}_{k+1}|^2 \, dx$$

$$= C_{(L)}^2 \|f - \bar{f}_{k_0}\|_{(1-\beta),2,2}^2.$$

Therefore, for $\mathcal{F}$ bounded in Lipschitz($\beta$), a compact Sobolev($\alpha$) set is also compact in $\mathcal{B}_{1/2,4,4}$ if $\alpha \geq 1 - \beta$.

Finally, if $\mathcal{F}$ is compact in Sobolev($\alpha$), $\alpha \geq 3/4$, then it immediately follows from the Sobolev embedding theorem that the function is bounded in Lipschitz(1/4) [e.g., Folland (1984), pages 270 and 273], and it follows that $\mathcal{F}$ is compact in $\mathcal{B}_{1/2,4,4}$.  $\square$

**Acknowledgments.** We thank the referees and an Associate Editor for several suggestions which led to improvements in the final manuscript.

L. D. Brown
M. G. Low
Department of Statistics
The Wharton School
University of Pennsylvania
Philadelphia, Pennsylvania 19104-6340
USA
e-mail: lbrown@wharton.upenn.edu
e-mail: lowm@wharton.upenn.edu

A. V. Carter
Department of Statistics
  and Applied Probability
University of California, Santa Barbara
Santa Barbara, California 93106-3110
USA
e-mail: carter@pstat.ucsb.edu

C.-H. Zhang
Department of Statistics
504 Hill Center, Busch Campus
Rutgers University
Piscataway, New Jersey 08854-8019
USA